\documentclass[numreferences]{kluwer}
\usepackage{graphicx}
\usepackage{setspace}
\usepackage{enumitem}
\usepackage{amsthm}
\newtheorem{Thm}{Theorem}[section]

\newtheorem{Alg}[Thm]{Algorithm}

\begin{document}
	%\begin{article}
	\begin{opening}
		\title{Trigonometric Interpolation Based Approach for Second Order Fredholm Integro-Differential Equations}
		\author{Xiaorong \surname{Zou}\email{xiaorong.zou@bofa.com}}
		%\author{Xiaorong \surname{Zou}\email{xiaorzou@gmail.com}}
		\institute{Global Market Risk Analytic, Bank of America}
		\runningauthor{X. Zou}
		\runningtitle{TIBA for Second Order FIDE}
		\date{May 6, 2025}
		\classification{MSC2000}{Primary 65T40; Secondary 45B05}

\begin{abstract}
A trigonometric interpolation algorithm for non-periodic functions has been recently proposed and applied to study general ordinary differential equation (ODE).  This paper enhances the algorithm to approximate functions in $2$ dim space. Performance of the enhanced algorithm is expected to be similar as in $1$-dim case and achieve accuracy aligned with smoothness of the target function, which is confirmed by numerical examples.

As an application, the $2$-dim trigonometric interpolation method is used to develop an algorithm for the solution of a second order Fredholm integro-differential equation (FIDE).  There are several advantages of the algorithm. First of all, it converges quickly and high accuracy can be achieved with moderate size of grid points; Secondly, it can effectively address singularities of kernel functions and work well with general boundary conditions.  Finally, it can be enhanced to copy with other IDE such as Volterra IDE or IDE with high order ODE component.  The tests conducted in this paper include various boundary conditions with both continuous kernels and integrable ones with singularity.  Decent performance is observed across all covered scenarios with moderate size of grid points.			
\end{abstract}
\end{opening}
\section{Introduction}\label{sec:intro}
A new trigonometric interpolation algorithm was recently introduced in \cite{zou_tri}.  It leverages Fast Fourier Transformation (FFT) to achieve optimal computational efficiency and converges at speed aligned with smoothness of underlying function. In addition, it can be used to approximate non-periodic functions defined on bounded intervals.  Considering the analytic attractiveness of trigonometric polynomial,  especially in handling differential and integral operations, the proposed trigonometric estimation of a general function is expected to be used in a wide spectrum. In this paper, we continue on applications of the trigonometric interpolation algorithm to solve integro-differential equations (IDEs) that contains an integral and derivative of unknown target function. 

IDEs appear in mathematical physics, applied mathematics and engineering (\cite{ide_23}-\cite{ide_36}). In general, analytical solutions of IDEs are not feasible and certain approximations are required. 
The concepts of IDEs have motivated a large size of research work for numerical solutions in recent years. 
%Quite rich research works have been done 
Several methods have been established and
among them are the Haar Wavelet method \cite{ide_28},  the Tau method \cite{ide_10}, the Taylor polynomial (\cite{ide_38}), the Spline collocation method \cite{ide_1} and Chebyshev polynomial method (\cite{ide_4}, \cite{ide_6}, \cite{ode_desa}).

Fredholm IDE (FIDE) and Volterra IDE (VIDE) are two most popular types of IDEs. The difference is that integration region of FIDE is fixed while VIDE bears a variable integration region. We shall address VIDE separately in \cite{zou_tri_VIDE} and aim to solve the following FIDE in this paper.  
\begin{equation}\label{eq:linear_ode_order2}  %{eq:bdp_orde2}
	y''(x)=p(x)y'+q(x)y(x) +r(x) + \mu(x)\int^e_sk(x,t)y(t)dt, \quad x\in [s,e].  
\end{equation}
\begin{eqnarray}
	d_{11}y(s) + d_{12}y'(s) +d_{13}y(e) + d_{14}y'(e) &=& \alpha , \label{eq:nonlinear_ode_order2:diri}\\ 
	d_{21}y(s) + d_{22}y'(s) +d_{23}y(e) + d_{24}y'(e) &=& \beta, \label{eq:nonlinear_ode_order2:neum} 
\end{eqnarray}
where $p,q,r,\mu$ is continuous differential on $[s, e]$, the kernel function $k(x,y)$ is integrable with respect to $y$ on $[s,t]$ for any $x\in [s,e]$, the rank of matrix $D:=(d_{ij})_{1\le i\le 2, 1\le j\le 4}$ is $2$, and $\alpha, \beta$ are two real numbers. 

Recently, trigonometric interpolation has been used to solve second order ODE with mixed boundary conditions \cite{zou_tri_III} and \cite{zou_tri_IV}. The proposed methods leverages the attractiveness of trigonometric polynomial for the connection among target function $y$ and its derivatives $y',y''$.  For non-linear ODE,  The coefficients of trigonometric polynomial can be identified effectively through optimization process whose gradient vectors can be carried out by fast Fourier transformation \cite{zou_tri_III}.  For linear ODE,  ODE can be converted to linear algebraic system and the coefficients can be solved directly. The method is flexible to handle general mixed boundary conditions and can be extended to high-order ODE \cite{zou_tri_IV}. 

We shall follow the same idea to solve FIDE (\ref{eq:linear_ode_order2}-\ref{eq:nonlinear_ode_order2:neum}). In order to address the integral component in Eq. (\ref{eq:linear_ode_order2}), we first enhance interpolation algorithm developed in \cite{zou_tri} for $2$-variable functions so a kernel $k(x,y)$ can be approximated by a $2$-dim trigonometrical polynomial. The algorithm developed in \cite{zou_tri_IV} can then be enhanced to find a trigonometric approximation of FIDE (\ref{eq:linear_ode_order2}-\ref{eq:nonlinear_ode_order2:neum}).

There are some advantage of proposed trigonometric interpolation method. It produces accurate approximation as shown by numerical examples conducted in Section \ref{sec:numerical_examples}; it can effectively address singularity issue in kernel function as described in Subsection \ref{sec:trig}; it works for various boundary conditions.  Finally, it can be enhanced to work for other IDE such as Volterra IDE or IDEs with high order ODE component. 

The rest of paper is organized as follows.  In Section \ref{sec:trig}, we start with a brief review on  trigonometric interpolation algorithm developed in \cite{zou_tri}; then enhance it for two variable functions so it can be used to approximate kernel function $k(x,t)$ in Eq. (\ref{eq:linear_ode_order2}); the section is closed by numerical testing results on some kernels to demonstrate performance of the enhanced algorithm. Section \ref{sec:fide} is denoted to develop Algorithm \ref{Alg:FIDE} to convert FIDE (\ref{eq:linear_ode_order2}-\ref{eq:nonlinear_ode_order2:neum}) into a linear algebraic system. It begins  with reviewing related results in \cite{zou_tri_IV} to attack the ODE component in FIDE.  Integral component is addressed in Subsection \ref{subsec:fide} and \ref{subsec:fide_en}  for continuous and integrable kernels respectively.  The proposed algorithm is summarized in Subsection \ref{subsec:transform_linear_algebraic_system}. Section \ref{sec:numerical_examples} demonstrates quite a few numerical examples to assess the performance of Algorithm \ref{Alg:FIDE}. Tests covers four types of boundary conditions with both continuous kernels and integrable kernels with singularity. The results show decent accuracy for all covered scenarios with moderate grid point size. The summary is made on Section \ref{sec:summary}.

\section{Trigonometric Interpolation on Non-Periodic Functions}\label{sec:trig}
\subsection{Trigonometric Interpolation with one-variable Functions}\label{sec:trig_dim1}
In this subsection, we review relevant results of trigonometric interpolation algorithm developed in \cite{zou_tri} starting with following interpolation algorithm on periodic functions. 
\begin{Alg}\label{main_thm} Let $f(x)$ be an odd periodic function \footnote{Similar results for even periodic function is also available in \cite{zou_tri}.} with period $2b$ and $N=2M=2^{q+1}$ for some integer $q\ge 1$ and $x_j, y_j$ are defined by	
	\begin{eqnarray}
		x_j &:=& -b + j\lambda,  \quad \lambda= \frac{2b}N , \quad 0\le j <N, \label{x_grid_N} \\
		y_j &:=& f(x_j).  \label{f_N_interpolation_new} 
	\end{eqnarray}
	Then there is a unique $M-1$ degree trigonometric polynomial
	\begin{eqnarray*}
		f_M(x) &=& \sum_{0< j <M}a_j \sin\frac{j\pi x}b, \label{f_M_interpolation_new_odd}\\
		a_j&=&\frac{2}N\sum_{0< k <N} (-1)^j y_k \sin\frac{2\pi j k}{N}, \quad 0< j <M \label{aj_odd}
	\end{eqnarray*}
	such that it fits to all grid points, i.e.
	\[%\label{error_even_new_odd} 
		f_M(x_{k})=y_{k}, \quad  0\le k <N.
	\]
\end{Alg} 
One can computer coefficients by Inverse Fast Fourier Transform (ifft):
\[	
\{a_j (-1)^j\}_0^{N-1} = 2\times Imag (ifft(\{y_k\}_{k=0}^{N-1})).
\]
%Note that the uniqueness of solution provides the theoretic support for the algorithm proposed in Section \ref{sec:algorithm}.
%It is shown in \cite{zou_tri_I} that  $f_M(x)$ and associated derivatives converge in a desired order. 
%In \cite{zou_tri_II}, 
Algorithm \ref{main_thm} has been enhanced so it can be applied to a nonperiodic function $f$ whose $K+1$-th derivative $f^{(K+1)}(x)$ exists over a bounded interval $[s,e]$. To seek for a periodic extension with same smoothness, we assume that $f$ can be extended smoothly such that $f^{(K+1)}$ exists and is bounded over $[s-\delta, e+\delta]$ for certain $\delta>0$.  Such a periodic extension of $f$ can be achieved by a cut-off smooth function $h(x)$ with following property: 
%One can use the following method to periodically extend $g$ by a cut-off function $h$ with the following properties 
\begin{eqnarray*}
	h(x)=\left\{\begin{array}{cc}
		1 & x\in [s,e], \\
		0 & \mbox{$x<s-\delta$ or $x>e+\delta$}. \\
	\end{array}\right.
\end{eqnarray*}
A cut-off function with closed-form analytic expression is proposed  in \cite{zou_tri}.  Let
\begin{equation}\label{ob}
	o=s-\delta, \quad b=e+\delta -o,
\end{equation}
and define $F(x):=h(x+o)f(x+o)$ for $x\in [0,b]$. One can treat $F(x)$ as an odd periodic function with period $2b$. Apply Algorithm \ref{main_thm} to generate the trigonometric interpolation of degree $M-1$ with $N$ evenly-spaced grid points over $[-b,b]$
\[
F_M(x) = \sum_{0 < j < M} a_j \sin\frac{j\pi x}{b},
\]
and let 
\[%\label{fMext}
	\hat{f}_{M}(x)=F_M(x-o)=\sum_{0< j < M} a_j \sin\frac{j\pi (x-o)}{b}. 
%\end{equation*}
\]
$\hat{f}_{M}(x)|_{[s,e]}$ can be treated as an trigonometric interpolation of $f$ since $\hat{f}_{M}(x_k)=f(x_k)$ for all grid points $x_k\in [s,e]$.  Numerical tests on certain functions demonstrate that $\hat{f}$ approaches to $f$ with decent accuracy \cite{zou_tri}. 

\subsection{Trigonometric Interpolation with two-variables Functions}\label{sec:trig_dim2}
Let $z=f(x,y)$ be a continuous function on $[s_1-\delta_1,e_1+\delta_1]\times [s_2-\delta_2,e_2+\delta_2]$ and define
\[
o_1=s_1-\delta_1, \quad b_1=e_1+\delta_1 -o_1, \qquad o_2=s_2-\delta_2, \quad b_2=e_2+\delta_2 -o_2. 
\]
As in one dimension case, one can extend $f(x,y)$ to a periodic function.
\begin{Alg}\label{alg:dim2} Select two integers $\{q_i \ge 1\}_{i=1,2}$ and let $N_i=2M_i=2^{q_i+1}$.
	\begin{enumerate}
		\item Define the evenly-spaced grid mesh $x_{j}, y_{k}, z_{jk}$ as follows
		\begin{eqnarray*}
			x_{j} &:=& -b_1 + j\lambda,  \quad \lambda= \frac{2b_1}{N_1} , \quad 0\le j <N_1,  \\
			y_{k} &:=& -b_2 + k\lambda,  \quad \lambda= \frac{2b_2}{N_2} , \quad 0\le k <N_2, \\
			z_{jk} &:=& f(x_{j},y_{k}),  
		\end{eqnarray*}
		
\item Construct cut-off functions $h_1(x)$ and $h_2(y)$ over $[-b_1,b_1]\times [-b_2,b_2]$ respectively:
\begin{eqnarray*}
	h_1(x)=\left\{\begin{array}{cc}
		1 & x\in [s_1,e_1], \\
		0 & \mbox{otherwise}. \\
	\end{array}\right. \quad
	h_2(y)=\left\{\begin{array}{cc}
		1 & y\in [s_2,e_2], \\
		0 & \mbox{otherwise}. \\
	\end{array}\right.
\end{eqnarray*}
\item Let $F(x,y)=h_1(x+o_1)h_2(y+o_2)f(x+o_1,y+o_2)$ and extend it to periodic function over $[-b_1,b_1]\times [-b_2,b_2]$.
\item Construct interpolation coefficients  $C_1(M_1,M_2)$ in $x$ direction. Let $ \vec{x}=(x_0,\dots, x_{N_1-1})$ and Let $C_1$ be the $M_1\times M_2$ matrix such that, for each $k\in [0, M_2-1]$, $C_1(:,k)$  is the interpolation coefficients determined by grid points $F(\vec{x},y_{k})$ based on the method described in Section \ref{sec:trig_dim1}. The computational cost of this step is $N_2 N_1\log_2(N_1)$. 
\item  Let $\vec{y}=(y_{0},\dots, y_{N_2-1})$. The desired $2$-dim coefficient matrix $C=(c_{jk})$ can be constructed such that, for each $j\in [0,M_1-1]$, its $j$-th row is equal to the interpolation coefficients determined by grid points determined by $\vec{y}$ and  $C_1(j,:)$.  The computational cost of this step is $N_1 N_2\log_2(N_2)$. 
\item Construct the following function that interpolates $f(x,y)$ at grid points
\[
\hat{f}_{M_1,M_2}(x,y) = \sum_{0<j<M_1,0<k<M_2}c_{jk}\sin\frac{j\pi (x-o_1)}{b_1}\sin\frac{k\pi (y-o_2)}{b_2}
\]
\end{enumerate} 
\end{Alg}
Note that the cost of Algorithm \ref{alg:dim2} is $N_1N_2\log_2(N_1N_2)$. 

\begin{figure}[htbp]\label{fig:alg_dim2}
	\centering
	\includegraphics[width=5cm]{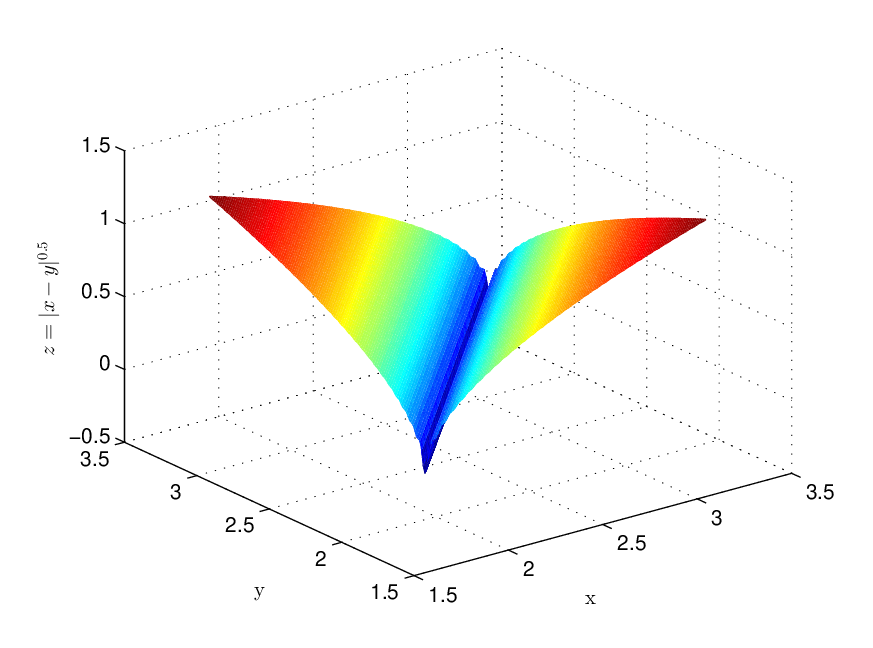}
	\includegraphics[width=5cm]{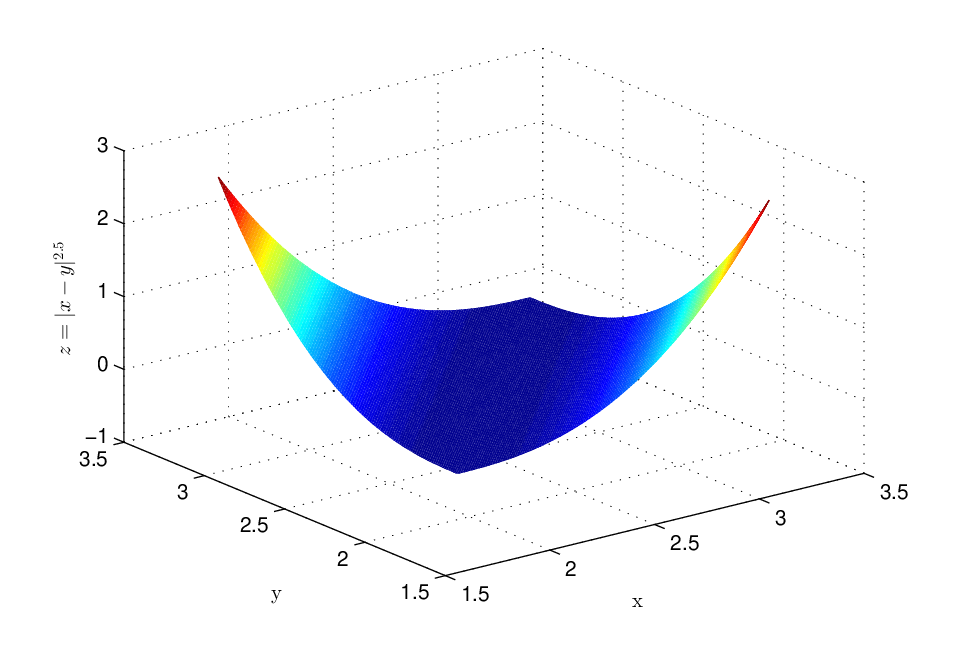}
	\caption{The plot of trigonometric interpolants for $|x-y|^{\gamma}$ with $\gamma=0.5$ (left) and $\gamma=2.5$ (right). Note that the figures are constructed by $2^{10}\times 2^{10}$ evenly-spaced points while interpolation is based on $M_1=M_2=2^7$.}
\end{figure}

The convergence performance of Algorithm \ref{alg:dim2} depends on smoothness of underlying function $f(x,y)$ and the number of grid points as in one-dimension situation. Tests are conducted on four types of functions as shown in Table \ref{tab:labels_test_dim2} based on the following setting 
\[
s_1=s_2=2, \quad  e_1=e_2=3, \quad  \delta_1=\delta_2=1, \quad b_1=b_2=3, \quad o_1=o_2 = 1.
\]
For type $K_1$ and $K_2$, the target function $k(x)$ is continuous, but not differentiable along the line $y=x$ for $\gamma<1$. The smoothness gets improved when $\gamma$ increases and is larger than $1$ as shown in Figure \ref{fig:alg_dim2}.  
\begin{table}[htbp]
\caption{The selected functions for the performance test on Algorithm \ref{alg:dim2}}
\begin{tabular}{cccc}
	Label & $k(x,y)$ & Label & 	$k(x,y)$ \\ \hline\hline
	$K_1$  & $|x-y|^\gamma$ & $Exp$  &  $exp(x+y)$  \\\hline
	$K_2$  & $|x^2-y^2|^\gamma$ & $Sin$ &  $sin(x+y)$ \\\hline  
\end{tabular}%
\label{tab:labels_test_dim2}%
\end{table}%
Performance of Algorithm \ref{alg:dim2} is measured by two normalized max errors over grid points and an extended set including non-grid points defined by
\begin{eqnarray*}
err_g &=& \frac{\max_{{x,y}\in S_g}|k(x,y)-\hat{k}_{M,M}(x,y)|}{\max_{(x,y)\in S_g}|k(x,y)|}, \\  
S_g &=& \{ (x_j, y_k), 0\le j,k \le 2^{q}  \} \cap [s,e]^2\\
err_e &=& \frac{\max_{{x,y}\in S_e}|k(x,y)-\hat{k}_{M,M}(x,y)|}{\max_{(x,y)\in S_e}|k(x,y)|}, \\ 
S_e &=& \{ (x^p_j, y^p_k), 0\le j,k \le 2^{10}  \} \cap [s,e]^2,\\
\end{eqnarray*}
where $(x^p_j,y^p_k)$ is equally spaced points based on $q_1=q_2=10$. 
Note that $S_g$ is part of the grid point set used in interpolation algorithm and is a subset of $S_e$ since $q<10$ in all test cases. 
$K_1$ is used as a testing kernel function in Section \ref{sec:numerical_examples} and Table \ref{tab:two-dim-convergence_k1} reports the two metrics with various $\gamma$ and $q$.  One can observe 
\begin{enumerate}
	\item $err_g$ is negligible under all cases as expected since the interpolant should match function values at $S_g$. 
	\item For $\gamma=0.5$, due to lack of smoothness,  $\hat{k}_{M,M}$ converges quite slowly and error can reach to $2.3E-02$ with $q=9$.  
	\item The performance  of  $\hat{k}_{M,M}$ gets improved significantly when $\gamma$ increases.  Decent accuracy is observed at $q=7$ for $\gamma=1.5$, which is supportive for the enhanced method developed in Section \ref{sec:trig} to handle kernels with singularities like $K_1$ and $K_2$.  See Section \ref{sec:numerical_examples} for numerical examples. 
\end{enumerate}   
\begin{table}[htbp]
	\caption{ Impact on Convergence with $K_1$}
	\begin{tabular}{lrrrr}
		case  & \multicolumn{1}{l}{$\gamma$} & \multicolumn{1}{l}{${\quad err_g}$} & \multicolumn{1}{l}{${\quad err_e}$} &  \multicolumn{1}{l}{\quad q} \\ \hline\hline
		K1    & 0.5   & 1.3E-15 & 6.6E-02 & 6 \\ \hline
		K1    & 0.5   & 2.0E-15 & 4.6E-02 & 7 \\ \hline
		K1    & 0.5   & 3.5E-15 & 3.2E-02 & 8 \\ \hline
		K1    & 0.5   & 3.3E-15 & 2.3E-02 & 9 \\ \hline
		K1    & 1.5   & 1.0E-15 & 3.8E-04 & 6 \\ \hline
		K1    & 1.5   & 1.9E-15 & 1.3E-04 & 7 \\ \hline
		K1    & 1.5   & 2.2E-15 & 4.6E-05 & 8 \\ \hline
		K1    & 1.5   & 2.3E-15 & 1.6E-05 & 9 \\ \hline
		K1    & 2.5   & 1.3E-15 & 2.9E-05 & 6 \\ \hline
		K1    & 2.5   & 1.5E-15 & 9.6E-07 & 7 \\ \hline
		K1    & 2.5   & 2.3E-15 & 1.7E-07 & 8 \\ \hline
		K1    & 2.5   & 2.4E-15 & 2.9E-08 & 9 \\ \hline
	\end{tabular}%
	\label{tab:two-dim-convergence_k1}%
\end{table}%

Table \ref{tab:K-1-K-2} compares the performance between $K_1$ and $K_2$ and similar results are observed.  As such,  we expect that similar performance of solving FIDE with kernel $K_1$ in Section \ref{sec:numerical_examples} should be applied to FIDE with kernel $K_2$ that is also popular in $IDEs$. 
\begin{table}[htbp]
	\caption{Performance Comparison between $K_1$ and $K_2$}
	\begin{tabular}{lrrrr}
		case  & \multicolumn{1}{l}{gamma} & \multicolumn{1}{l}{$err_g$} & \multicolumn{1}{l}{$err_e$} & \multicolumn{1}{l}{q} \\ \hline\hline
		K1    & 0.5   & 2.0E-15 & 4.6E-02 & 7 \\ \hline
		K2    & 0.5   & 2.6E-15 & 5.2E-02 & 7 \\\hline
		K1    & 1.5   & 1.9E-15 & 1.3E-04 & 7 \\\hline
		K2    & 1.5   & 1.8E-15 & 1.9E-04 & 7 \\\hline
		K1    & 2.5   & 1.5E-15 & 9.6E-07 & 7 \\\hline
		K2    & 2.5   & 1.7E-15 & 1.8E-06 & 7 \\\hline
	\end{tabular}%
	\label{tab:K-1-K-2}%
\end{table}%
Table \ref{tab:two-dim-four-functions} compares performance of four tested kernel function with $q=7$. Decent accuracy are observed for two smooth kernel $Exp$ and $Sin$ as expected. 
\begin{table}[htbp]
	\caption{Performance of four test functions with $q=7, \gamma=0.5$}
	\begin{tabular}{lrrrr}
		case  & \multicolumn{1}{l}{$\gamma$} & \multicolumn{1}{l}{{$err_g$}} & \multicolumn{1}{l}{$err_e$} & \multicolumn{1}{l}{q} \\ \hline\hline
		K1    & 0.5   & 2.0E-15 & 4.6E-02 & 7 \\\hline
		K2    & 0.5   & 2.6E-15 & 5.2E-02 & 7 \\\hline
		exp   & NA   & 1.7E-15 & 5.9E-08 & 7 \\\hline
		sin   & NA   & 2.0E-15 & 5.4E-08 & 7 \\\hline
	\end{tabular}%
	\label{tab:two-dim-four-functions}%
\end{table}%

\section{Trigonometric Approximation for Solutions of FIDF}\label{sec:fide}
Let $p(x)$,$q(x)$, $r(x)$, $K(x,t)$ be defined in Eq (\ref{eq:linear_ode_order2}), define
\[
f(x,v,u) = p(x)u+q(x)v +r(x), \quad  g(x,v) = \mu(x)\int^e_sK(x,t)v(t)dt
\]
FIDE (\ref{eq:linear_ode_order2}) is reduced to linear ODE if the integration component $g(x,v)$ vanishes and has been solved in \cite{zou_tri_III} by transforming ODE and boundary condition (\ref{eq:nonlinear_ode_order2:diri}-\ref{eq:nonlinear_ode_order2:neum}) into linear algebraic system.  We shall also convert FIDE (\ref{eq:linear_ode_order2})  to a linear algebraic system by leveraging what has been developed in \cite{zou_tri_III} to handle  ODE component $f$ and attack integration component $g$ in Subsection \ref{subsec:fide} and \ref{subsec:fide_en} depending on smoothness of kernel function $k$.

To apply trigonometric interpolation for non-periodic functions, 
we assume that $p(x)$, $q(x)$ and $r(x)$ are continuous differential on $[s-\delta, e+\delta]$ for certain $\delta>0$ and $K(x,t)$ are continuous differential on $[s-\delta, e+\delta]^2$. 
By parallel shifting if needed,  we assume $s=\delta$ without loss of generality. Let $h$  be a cut-off function specified in Section \ref{sec:trig} and 
%extend $g(x),p(x),q(x)$ to $g_h(x),p_h(x), q_h(x)$ as follows
construct $f_h(x,v,u)$ and $g_h(x,t)$ as follows
\[
f_h(x,v,u) = f(x,v,u)h(x), \quad g_h(x,v) =h(x)g(x,v).
\]
Consider a solution $v(x)$ of the following FIDE system
\begin{eqnarray}
	%v''(x) &=& g_h(x) +p_h(x) v'(x) + q_h(x) v(x),  \quad x\in [0,b] \label{eq:nonlinear_ode_order2_F},\\
	v''(x) &=& f_h(x, v,v') + g_h(x,v),  \quad x\in [0,b] \label{eq:nonlinear_ode_order2_F},\\
	\alpha  &=& d_{11}v(s) + d_{12}v'(s) +d_{13}v(e) + d_{14}v'(e), \label{eq:nonlinear_ode_order2:diri_F}\\ % Dirichlet boundary condition 
	\beta  &=& d_{21}v(s) + d_{22}v'(s) +d_{23}v(e) + d_{24}v'(e). \label{eq:nonlinear_ode_order2:neum_F} 
\end{eqnarray}
It is clear that $v(x)|_{[s,e]}$ solves FIDE (\ref{eq:linear_ode_order2}-\ref{eq:nonlinear_ode_order2:neum}). Define $u(x):=v'(x)$ and $z(x):=v''(x)$.  By Eq (\ref{eq:nonlinear_ode_order2_F}),  $z(x)$ and its derivatives $z^{(k)}$ vanish at boundary points $\{0,b \}$, hence it can be smoothly extended as an odd periodic function with period $2b$ and be approximated by trigonometric polynomial.  

Assume 
\begin{equation}\label{eq:z_M_nonlinear}
	\tilde{z}_M(x) = \sum_{0\le j<M}b_j \sin \frac{j\pi x}{b} 
\end{equation}
is an interpolant of $z(x)$ with $N$ equispaced grid points over $[-b,b]$. $u$ and $v$ can be derived accordingly 
\begin{eqnarray}
	\tilde{u}_M(x ) &=& a_0 -   \frac{b}{\pi }\sum_{1\le j <M} \frac{b_j}{j} \cos  \frac{j\pi x}{b},  \label{eq:tildeu}\\
	\tilde{v}_M(x ) &=& a_{-1} +  a_0 x -  (\frac{b}{\pi })^2 \sum_{1\le j <M} \frac{b_j}{j^2}  \sin  \frac{j\pi x}{b},  \label{eq:tildev_nonlinear}
\end{eqnarray}
where $a_0,a_{-1}$ are two constant and can be determined by boundary conditions as shown in Eq (\ref{eq:a01}) below. 

The following notations and conventions will be adopted in the rest of this paper. A $k$-dim vector is considered as $(k,1)$ dimensional matrix unless specified otherwise. Define
\begin{eqnarray*} 
	x_k &=& k \lambda, \quad \lambda = \frac{b}M ,\quad X = (x_k)_{1\le k < M},  \\
	u_k &=& \tilde{u}_M(x_k), \quad v_k = \tilde{v}_M(x_k),  \quad z_k =\tilde{z}_M(x_k), \quad \mu_k = \mu(x_k)h(x_k),\\
	p_k &=& p(x_k)h(x_k), \quad q_k=q(x_k)h(x_k), \quad r_k = r(x_k)h(x_k), \\
	f_k &=& f_h(x_k, v_k, u_k),  \quad g_k =g_h(x_k, v(s:e)), \\
	%g_k &=& g_h(x_k), \quad p_{k} = p_h(x_k), \quad q_k = q_h(x_k), \quad r_k= g_k+p_ku_k+q_kv_k\\
	U &=& (u_k),\quad V = (v_k), \quad Z = (z_k),  \quad 1\le k< M,\\
	F &=& (f_k), \quad G = (g_k) , \quad B=(b_i), \quad 1\le k < M\\
	I &=& (1,1,\cdots)^T_{M-1}, \quad I_a = (-1,1,\cdots)^T_{M-1}, \quad K = (1,2,\cdots, M-1)^T, \\
	P &=& (p_k), \quad Q = (q_k), \quad R= (r_k), \quad {\nu} = (\mu_k) \quad 1\le k< M.
\end{eqnarray*}
For any two matrices $A,B$ with same shape,  $A\circ B$ denotes the Hadamard product, which applies the element-wise product to two matrices. $A\cdot B$ denote the standard matrix multiplication when applicable. $A(i,:)$ and $A(:,j)$ is used to denote the $i$-th row and $j$-th column of $A$ respectively.  $W(k:l)$ denote the vector $(w_k, \cdots, w_l)^T$. In addition, $diag (W)$ is the diagonal matrix constructed by $W$.  Note we have
\[
s = x_{m}, \qquad e=x_{m+n}.
\]
At grid points of interpolation, ODE dynamic (\ref{eq:nonlinear_ode_order2_F}) is characterized by
\begin{equation} \label{eq:bdp_orde2_d_nonlinear} 
	Z = F + G.  
\end{equation}
$\{z_k, u_k, v_k\}_{0\le k \le M}$ can be calculated based on Eq. (\ref{eq:z_M_nonlinear}-\ref{eq:tildev_nonlinear}):
\begin{eqnarray}
	z_k &=& \sum_{0\le j<M}b_j \sin \frac{2\pi jk}{N},   \label{eq:z_M_d}\\
	u_k &=& a_0 - \frac{b}{\pi }\sum_{1\le j <M} \frac{b_j}{j} \cos  \frac{2\pi jk }{N},  \label{eq:tildeu_d}\\
	v_k &=& a_{-1} + a_0 x_k-(\frac{b}{\pi })^2 \sum_{1\le j <M}\frac{b_j}{j^2}\sin \frac{2\pi jk }{N}. \label{eq:tildev_d}
\end{eqnarray}
Eq (\ref{eq:tildev_d}) can be used to solve $a_0$ and $a_{-1}$:
\begin{equation}\label{eq:a01}
	a_0 = \frac{v_M - v_0}b, \quad a_{-1} = v_0.
\end{equation}
Define 
%\begin{equation*}\label{eq:SC}
\[
	S =(\sin\frac{2\pi jk}{N})_{1\le j,k <M}, \qquad C =  (\cos\frac{2\pi jk}{N})_{1\le j,k < M},
\]
%\end{equation*}
and 
\[ 
 O=\sqrt{\frac M2}S,\qquad \Theta = O \cdot diag(1/K^2), \quad  O := (\theta_{ij})_{1\le i, j < M}.
\]
$B$ can be solved by $V$ as
\begin{equation}\label{eq:B}
	B = diag(K^2) S (\frac{2a_{-1} \pi^2}{Mb^2} I + \frac{2a_0 \pi^2}{bM^2}K - \frac{2\pi^2}{Mb^2} V).
\end{equation}
Let $U_e=(u_0,\dots, u_M)^T$, $V_e=(v_0,\dots, v_M)^T$, it is shown in  \cite{zou_tri_IV} that $U_e$ is covered by $V_e$ through linear transform   
\begin{equation}\label{eq:UAV}
	U_e = AV_e,
\end{equation}
where matrix $A$ is calculated by
\begin{eqnarray*}
	a_{0,0} &=& \frac{\pi}{b} sum(I_a \circ \cot(\pi K/N))\\
	&-&\frac{\pi}{bM} sum(I_a\circ K\circ \cot(\pi K/N))-\frac1b \label{eq:top}\\
	a_{0,1:M-1} &=& -\frac{\pi}{b} I'_a \circ \cot(\pi K'/N) \nonumber \\
	a_{0,M} &=& \frac{\pi}{bM} sum(I_a\circ K\circ \cot(\pi K/N))+\frac1b \nonumber\\
	a_{i,0} &=& \frac{\pi}{2b} sum((-1)^i \cot(i,:)I_a) \\
	&-& \frac{\pi}{2bM} sum( (-1)^i I_a \circ \cot(i,:) \circ K) -1/b, \label{eq:middle}\\
	a_{i,1:M-1} &=& \frac{\pi}{2b}(-1)^{i+1} I_a' \circ \cot(i,:), \nonumber \\
	a_{i,M} &=& \frac{\pi}{2bM} sum((-1)^i I_a \circ \cot(i,:) \circ K)+1/b. \nonumber \\
	a_{M,0} &=& -\frac{\pi}{b} sum(I_a \circ \tan(\pi K/N))\\
	&-&\frac{\pi}{bM} sum((K\circ \cot(\pi K/N))-\frac1b, \label{eq:bottom}\\
	a_{M,1:M-1} &=& \frac{\pi}{b} I'_a \circ \tan(\pi K'/N), \nonumber\\
	a_{M,M} &=& \frac{\pi}{bM} sum((K\circ \cot(\pi K/N))+\frac1b,  \nonumber
\end{eqnarray*}
for $\quad 0<i<M$ and
\[
\cot(k, i) := Cot\frac{k+i}{N}\pi + Cot\frac{k-i}{N}\pi,
\]
and $Cot(x)=\cot(x)$ if $x/\pi$ is not integer and $Cot(x)=0 $ otherwise. Eq (\ref{eq:bdp_orde2_d_nonlinear}) is equivalent to 
\begin{equation}\label{eq:ode_discrete_nonlinear}
	\frac{v_0\pi^2}{Mb^2} (MI-K) + \frac{v_M\pi^2}{M b^2} K -\frac{\pi^2}{b^2} V = \Theta \cdot (R + Q\circ V  + P\circ U +  G).
\end{equation}

With Eq (\ref{eq:UAV}), Eq. (\ref{eq:ode_discrete_nonlinear}) represents a linear algebraic system if $G$  can be also covered by $V$ through linear transformation, which is to be developed in Subsection \ref{subsec:fide} and \ref{subsec:fide_en}. 

\subsection{Integral component with continuous kernel}\label{subsec:fide} 
We interpolate $k(x,y)$ by $2$-dim $sin$ polynomial as shown in Section \ref{sec:trig_dim2}   
\begin{equation}
	K(x,t) = \sum_{1\le i,j<M} c_{ij} \sin \frac{\pi i x}b\sin \frac{\pi j t}b,
\end{equation}
and apply Eq (\ref{eq:tildev_nonlinear}) to estimate $g_k$ for $1\le k <M$
\begin{eqnarray*}
	\frac{g_k}{\mu_k} &=& \sum_{1\le i,j<M} c_{ij} \int^e_s  \sin \frac{2\pi i k}N \sin \frac{\pi j t}b (a_1 +a_0 -\frac{b^2}{\pi^2} \sum_{1\le l<M} \frac{b_l}{l^2} \sin\frac{\pi l t}{b})dt ,\\
	&=& \sum_{1\le i,j<M} \frac {bc_{ij}}{\pi j} \sin\frac{2\pi ik}{N} (w^j_{-1}a_{-1} + w^j_0 a_0) \\
	&-&  \frac {b^3}{2\pi^3}  \sum_{1\le i,j,l<M} c_{ij} \sin\frac{2\pi ik}{N} b_l w^j_{l},   
\end{eqnarray*}
where
\begin{eqnarray*}
	w^j_{-1} &=& \cos\frac{2\pi j m}{N} -  \cos\frac{2\pi j (m+n)}{N}, \label{f_w_1}\\ 
	w^j_0 &=& s\cos\frac{2\pi j m}{N} -  e\cos\frac{2\pi j (m+n)}{N} + \frac{b}{\pi j} (\sin\frac{2\pi j (m+n)}{N} \\
	&-&\sin\frac{2\pi j m}{N}  ) \label{f_w_0}, \\
	w_{jl} &=& \frac{\delta_{n\neq l}}{(j-l)l^2} (\sin\frac{2\pi (j-l) (m+n)}{N} -\sin\frac{2\pi (j-l) m}{N} ) \label{f_w}\\
	&+& \frac{2\pi n}{l^2 N}\delta_{j=l} -\frac{1}{(j+l)l^2} (\sin\frac{2\pi (j+l) (m+n)}{N} -\sin\frac{2\pi (j+l) m}{N} ).
\end{eqnarray*}
Define 
\[
w_{-1}=(w^j_{-1})_{1\le j<M}, \quad w_{0}=(w^j_{0})_{1\le j<M}, \quad  W=(w_{jl})_{1\le j, l<M}.
\]
Let 
\begin{eqnarray*}
	\hat \eta &=& (\hat \eta_{kj})_{(M-1, M-1)} , \quad \hat \eta_{kj} = \sum_{0<i<M}c_{ij}/j \sin\frac{2\pi ik}{N},  \\
	%\hat \eta_{j} &:=& (\hat \eta_{kj})_{1\le k<M}= N\times imag( ifft(c_{:,j}/j)(1:M-1),  \\ 
	\eta &=& (\eta_{kj})_{(M-1, M-1)} , \quad \eta_{kj} = \sum_{0<i<M}c_{ij} \sin\frac{2\pi ik}{N}. 
	%\eta_{j} &:=& (\eta_{kj})_{1\le k<M} = N\times imag( ifft(c_{:,j}))(1:M-1),\\	
\end{eqnarray*}
We have
\begin{eqnarray*}
	\frac{g_k}{\mu_k}&=& \sum_{1 \le j<M} \frac {b \hat\eta_{kj}}{\pi }  (w^j_{-1}a_{-1} + w^j_0 a_0) -\sum_{1\le j,l<M}  \frac {b^3 \eta_{kj}}{2\pi^3} b_l w_{jl},   
\end{eqnarray*}
or
\begin{equation}\label{f-govermu}
	G/\nu = \frac b{\pi} \hat\eta w_{-1}a_{-1}  + \frac b{\pi} \hat\eta w_{0} a_0 - \frac{b^3}{2\pi^3} \eta WB. 
\end{equation}
Combining Eq (\ref{eq:a01}-\ref{eq:B}) and (\ref{f-govermu}),  $G$ can be recovered by $V$ as follows   
\begin{equation}\label{G_v_continuous}
	G=C_0 v_0 +C_M v_M  +  C_{rest}V,
\end{equation}
where
\begin{eqnarray*}
	H &=& \eta\cdot \dot W\cdot diag(K^2)\cdot S\\
	C_M &=& \frac{ b}{\pi} \nu\circ  ( \hat\eta \cdot W_0/b - H\cdot K/M^2) \nonumber\\
	C_0 &=& \frac{ b}{\pi} \nu\circ  (\hat\eta \cdot W_{-1} - H\cdot I/M)  -C_M \nonumber \\
	C_{rest} &=& \frac{ b}{\pi} \cdot  diag (\nu)\cdot H/M. \nonumber
\end{eqnarray*}
\subsection{Integral component with integrable kernel}\label{subsec:fide_en} 
We assume that there are some singularity such that $k(x,y)=(|x-y|^{\gamma})\kappa(x,y)$ for $\gamma>-1$ and there are $k_1(x,y)$ and $k_2(x,y)$ such that $k_1(x,x)=0$ and $k_2(x,x)=0$ and 
\[
\frac{\partial {k_1(x,y)}}{\partial y}= k(x,y), \quad \frac{\partial {k_2(x,y)}}{\partial y}= k_1(x,y)
\]
As example $k(x,y)=|x-y|^{\gamma}$, we have
\begin{eqnarray*}
	k_1(x,y)=\left\{\begin{array}{cc}
		-\frac{|x-y|^{1+\gamma}}{1+\gamma} & y\le x, \\
		\frac{|x-y|^{1+\gamma}}{1+\gamma} &  y> x\\
	\end{array}\right.  
\end{eqnarray*}
\[
 k_2(x,y) = \frac{1}{(1+\gamma)(2+\gamma)}|x-y|^{(1+\gamma)(2+\gamma)}.
\]
So we have
\begin{equation}
	g(x,v) = -v(s)k_1(x,s)+u(s)k_2(x,s) + \mu(x) \int^x_s k_2(x,t)z(t).
\end{equation}
We extend $k_2(x,y)$ to $[-b,b]^2$ by $2$-dim $sin$ polynomial    
\[
	k_2(x,t) = \sum_{1\le i,j<M} c_{ij} \sin \frac{\pi i x}b\sin \frac{\pi j t}b,
\]
and apply Eq (\ref{eq:tildev_nonlinear}) to estimate $g_k$ for $1\le k <M$: 
\begin{equation}\label{v_g_integrable}
	\frac{g_k}{\nu_k}= -k_1(x_k, s) v_m  + k_2(x_k,s)u_m + \frac {b}{2\pi} \sum_{1\le j,l<M} \eta(k,j) w^k(j,l)b_l,    
\end{equation}
\begin{eqnarray*}
	\eta_{kj} &=& \sum_{i}c_{ij} \sin\frac{2\pi ik}{N},  \\
	w^k(j,l)&=& \frac{\delta_{j\neq l}}{(j-l)} (\sin\frac{2\pi (j-l) k}{N} -\sin\frac{2\pi (j-l) m}{N} ) \\
	&+&\frac{\pi}{b}(x_k-s)  -\frac{1}{ (j+l)} (\sin\frac{ 2\pi(j+l) k}{N} -\sin\frac{2\pi (j+l) m}{N} ).
\end{eqnarray*}
Define 
\begin{eqnarray*}
	\eta &=& (\eta_{kj})_{0< k, j<M},  \\
	%W_{-1} &=& (w_{-1}(k,j))_{(0< k, j<M)}, \quad W_{0}=(w^j_{k,0})_{(0<j, k<M)}, \\
	W^k &=& (w^k_{j,l})_{0< j, l<M}, \quad  H(k,:)=\eta(k,:)\cdot W^k\cdot diag(K^2)S, \quad 1\le k<M.
\end{eqnarray*}
For any matrix $A$, let $A_s$ denote the column vector constructed by the summation of all columns of $A$.  By Eq. (\ref{eq:a01}-\ref{eq:B}), Eq (\ref{v_g_integrable}) is reduced to
\begin{equation}\label{G_v_integrable}
	G=  C_{0} v_0 + C_mv_m + C_M v_M +  C_{rest} V,
\end{equation}
where
\begin{eqnarray*}
	C_0 &=& \frac{\pi}{b} \nu \circ (H\cdot I/M - H\cdot K/M^2) +  A(m,0) \cdot (\nu \circ k_2(X,s)),\\
	C_M &=& \frac{\pi}{bM^2} \nu \circ (H\cdot K) + A(m, M) \cdot (\nu \circ k_2(X,s)),\\
	C_m &=& -\nu\circ k_1(X,s), \\
	C_{rest} &=& -\frac{\pi}{b M} \cdot  diag(\nu)\cdot H  + (\nu \circ k_2(X,s))\cdot A(m, 1:M-1).\nonumber
\end{eqnarray*}
\subsection{Algorithm of solving FIDE}
\label{subsec:transform_linear_algebraic_system}

With linear representation of $G$ either by Eq (\ref{G_v_continuous}) or Eq (\ref{G_v_integrable}), we can solve FIDE (\ref{eq:linear_ode_order2},\ref{eq:nonlinear_ode_order2:diri}-\ref{eq:nonlinear_ode_order2:neum}) by Algorithm \ref{Alg:FIDE}.
\begin{Alg}\label{Alg:FIDE} Trigonometric Interpolation Based Approach for FIDE. 
\begin{enumerate}
\item Construct the following $M+1$ dimensional linear algebraic system  that consists of Eq (\ref{eq:nonlinear_ode_order2:diri_F}) (first in the system), Eq (\ref{eq:ode_discrete_nonlinear}) (last in the system), and Eq (\ref{G_v_continuous}) or Eq (\ref{G_v_integrable}). 
\begin{equation}\label{eq:PhiVG}
	\Phi V_e = \Psi,  \quad \Phi = (\phi_{jk})_{0\le j,k\le M}, \Psi = (\psi_j)_{0\le j\le M}.
\end{equation}
\item  Calculate $\Psi=(\psi_i)_{0\le i\le M}$ by 
\begin{equation} \label{eq:psi}
	\psi_0=\alpha, \qquad \psi_M= \beta, \quad  \Psi(1:M-1)= -\Theta R.
\end{equation}
\item Calculate $\Phi=(\phi_{ij})_{0\le i,k \le M}$ in three steps
\begin{enumerate}
	\item Identify the first and last row of $\Phi$, which are associated to Eq  (\ref{eq:nonlinear_ode_order2:diri_F}) and Eq (\ref{eq:ode_discrete_nonlinear}) respectively  %Find the $\Phi$ without the impact of $C_m$ and $C_{m+n}$: 
	\begin{eqnarray*}
		\phi_{0,k} &=& d_{12}\cdot A(m,k) + d_{14}\cdot A(m+n,k) \\
		&+& \delta_{m,k}d_{11} + \delta_{m+n,k}d_{13}, \label{eq:phi_0} \\
		\phi_{M,k} &=& d_{22}\cdot A(m,k) + d_{24}\cdot A(m+n,k) \\
		&+& \delta_{m,k}d_{21} + \delta_{m+n,k}d_{23}. \label{eq:phi_M} 	 
	\end{eqnarray*}
	\item Identify other $M-1$ rows of $\Phi$. For $0<i<M$, 
	\begin{eqnarray*}
		\phi_{i,0} &=& -\frac{(M-i)\pi^2}{Mb^2} + (\theta(i,:)\circ P^T) \cdot A(1:M-1,0) \\
		&+& \theta(i,:) C_0,\label{eq:phi_i_0} \\
		\phi_{i,M} &=& -\frac{i \pi^2}{b^2 M} + (\theta(i,:) \circ P^T) \cdot A(1:M-1, M) \\
		&+& \theta(i,:)C_M,\label{eq:phi_i_M} \\
		\phi(i,1:M-1) &=& \frac{\pi^2}{b^2}(\delta_{ij})_{1\le j<M} +\theta(i,:)\circ Q^T  +\nonumber\\
		&+&(\theta(i,:) \circ P^T) \cdot A(1:M-1,1:M-1) \label{Phi}  \\
		&+& \theta(i,:)C_{rest}  \label{eq:phi_i_middle}
	\end{eqnarray*}
	\item  Add the impact of $C_m$ if Eq (\ref{G_v_integrable})  is used to estimate $G$ component.
	\begin{eqnarray*}
	\Phi(1:M-1,m) \leftarrow \Phi(1:M-1,m) + \theta\cdot c_m
	\end{eqnarray*}
\end{enumerate}
\item Solve Eq (\ref{eq:PhiVG}) with calculated $\Phi$ and $\Psi$ in previous steps.
\item Apply Eq (\ref{eq:a01}) and Eq (\ref{eq:B}) to compute coefficients in Eq (\ref{eq:tildev_nonlinear}) for the approximation $v_M$ on solution of FIDE (\ref{eq:linear_ode_order2},\ref{eq:nonlinear_ode_order2:diri}-\ref{eq:nonlinear_ode_order2:neum}).
\end{enumerate}
\end{Alg}
\newpage
	
\section{Numerical examples}\label{sec:numerical_examples}
%to handle singularity of kernel
In this section,  we  test the method proposed in Section \ref{sec:fide} with kernel functions $k(x,y)$ defined in Table \ref{tab:test} with focus on $K_1$. 
%$k(x,y)=|x-y|^{\gamma}$ ($\gamma>-1$).  
For any given smooth functions $p(x),q(x),\mu(x)$, matrix $(d_{ij})_{2\times 4}$ and $f(x)\in C^2([s, e])$, one can see that $f(x)$ solves the following IDEF
\begin{eqnarray*}
	y'' &=& p(x) y' +  q(x)y  + r(x) + \mu(x) \int^e_sk(x,t)y(t)dt\\
	\alpha  &=& d_{11}y(s) + d_{12}y'(s) + d_{13}y(e) + d_{14}y'(e) \\
	\beta  &=& d_{21}y(s) + d_{22}y'(s) + d_{23}y(e) + d_{24}y'(e) 
\end{eqnarray*}
where $r(x)$ and constant $\alpha,\beta$ are determined by
\begin{eqnarray}
 r(x) &=& f''(x)  - p(x)f'(x) - q(x)f(x) -  \mu(x)\int^e_s k(x,t)f(t)dt \\
 	\alpha  &=& d_{11}f(s) + d_{12}f'(s) + d_{13}f(e) + d_{14}f'(e) \\
 \beta  &=& d_{21}f(s) + d_{22}f'(s) + d_{23}f(e) + d_{24}f'(e). 
\end{eqnarray}
All tests are based on the following parameters 
\[
s=1,\quad e=3, \quad \delta=1, \quad p(x) \equiv  0.1, \quad q(x) \equiv 1, \quad \mu(x) \equiv 1,
\]
target function $f$ can be one of following list 
\begin{equation}\label{test_functions}
\cos (\frac{\pi x}2 ) , \quad  \cos (\frac{3\pi x}2) , \quad e^x, \quad x^2. 
\end{equation}
We consider four sets of boundary conditions in Table \ref{tab:test}.
\begin{table}[htbp]
	\small
	\caption{The types of boundary conditions. $\{d_{ij}\}_{1\le i,j\le 4}$ are parameters in Eq (\ref{eq:nonlinear_ode_order2:diri}-\ref{eq:nonlinear_ode_order2:neum}).}
	\begin{tabular}{lrrrrrrrrl}
		type	& \multicolumn{1}{l}{$d_{11}$} & \multicolumn{1}{l}{$d_{12}$} & \multicolumn{1}{l}{$d_{13}$} & \multicolumn{1}{l}{$d_{14}$} & \multicolumn{1}{l}{$d_{21}$} & \multicolumn{1}{l}{$d_{22}$} & \multicolumn{1}{l}{$d_{23}$} & \multicolumn{1}{l}{$d_{24}$} & condition on \\ \hline\hline
		$Neumann$ & 1     & 0     & 0     & 0     & 0     & 1     & 0     & 0     & $v_s,u_s$ \\
		$Dirichlet$ & 1     & 0     & 0     & 0     & 0     & 0     & 1     & 0     & $v_s, v_e$ \\
		$Mix_1$   & 1     & 0     & 0     & 0     & 0     & 0     & 0     & 1     & $v_s, u_e$ \\ 
		$Mix_2$ & 1     & 1     & 0     & 0     & 0     & 0     & 1     & 1     & $v_s+u_s, v_e+u_e$ \\ \hline
	\end{tabular}%
	\label{tab:test}%
\end{table}%
The performance is measured by the normalized max error $err_{e}$ defined as
\[
max_{e} = \frac{\max_{x\in S_e}|f(x)-\tilde{v}_M(x)|}{\max_{x\in S_e}|f(x)|}, \quad  S_e = \{ x^p_k, 0\le k \le 2^{11}  \} \cap [s,e]
\]  
where $S_e=\{x^p_k\}$ consists equally spaced points over $[-b,b]$ based on $q=10$ and contains more than interpolation grid points used to determine $\tilde v_M$ since $M$ is always less than $2^{10}$ in all test cases. 
\subsection{Numerical results on performance of continuous kernels}\label{subsec:numerical_results_continous}
The performance of four continuous kernels in Table \ref{tab:labels_test_dim2} are tested with target function $f(x)=\cos(3\pi x/2)$ and Dirichlet boundary condition,  and the results are shown in Table \ref{tab:Fredholm_en_convergenc}.  Almost same accuracy is observed across four kernel functions, which is a little surprising since trigonometric interpolation error for $K_1$ and $K_2$ can reach to $E-02$, much higher than the other two kernel as shown in Table \ref{tab:two-dim-four-functions}.  The testing result might imply that the accuracy of FIDE's solution is mainly determined by the performance of trigonometric interpolation of kernel at grid points, where four kernels exhibits similar performance as shown in Table \ref{tab:two-dim-four-functions}. 
\begin{table}[htbp]
	\caption{The performance of covered kernels with   $f(x)=\cos (\frac{3\pi x}2)$}
	\begin{tabular}{lrrrrr}
		kernel & type  & \multicolumn{1}{l}{$\gamma$} & \multicolumn{1}{l}{$max_{e}$} &  \multicolumn{1}{l}{q} \\ \hline\hline
		$K_1$ & Dirichlet & 0.5   & 5.0E-11  & 7 \\ \hline
		$K_2$ & Dirichlet & 0.5  & 5.0E-11 &  7 \\ \hline
		$Exp$ & Dirichlet & NA   & 5.0E-11  & 7 \\ \hline
		$Sin$ & Dirichlet & NA   & 5.0E-11  & 7 \\ \hline
	\end{tabular}%
	\label{tab:Fredholm_en_convergenc}%
\end{table}%

Table \ref{tab:test_on_boundary} shows the performance of different boundary conditions with target function $f(x)=\cos(3\pi x/2)$ for FIDE with kernel $K_1$. The performance is decent for all cases and reaches the best accuracy with Dirichlet condition. 
\begin{table}[htbp]
	\caption{The performance of covered boundary conditions with   $f(x)=\cos (\frac{3\pi x}2)$}
	\begin{tabular}{lrrrrr}
		kernel & type  & \multicolumn{1}{l}{$\gamma$} & \multicolumn{1}{l}{$max_{e}$} & \multicolumn{1}{l}{q} \\ \hline\hline
		$K_1$ & Neumann & 0.5    & 3.5E-08 &  7 \\ \hline
		$K_1$ & Dirichlet & 0.5  & 5.0E-11 &  7 \\ \hline
		$K_1$ & mix\_1 & 0.5    & 3.6E-09 &  7 \\ \hline
		$K_1$ & mix\_2 & 0.5    & 1.9E-08 &  7 \\ \hline 
	\end{tabular}%
	\label{tab:test_on_boundary}%
\end{table}%

Table \ref{tab:test_on_convergence} shows convergence property and high accuracy is achieved for $q\ge 7$ with Dirichlet condition. As expected,  the performance is improved as $q$ increases. 
\begin{table}[htbp]
	\caption{Convergence Test with $f(x)=\cos (\frac{3\pi x}2)$}
	\begin{tabular}{lrrrr}
		kernel & type  & \multicolumn{1}{l}{$\gamma$} & \multicolumn{1}{l}{$max_{e}$} & \multicolumn{1}{l}{q} \\ \hline\hline
		$K_1$ & Dirichlet & 0.5   & 4.8E-03 & 4 \\ \hline
		$K_1$ & Dirichlet & 0.5   & 9.3E-05 & 5 \\ \hline
		$K_1$ & Dirichlet & 0.5   & 1.6E-07 & 6 \\ \hline
		$K_1$ & Dirichlet & 0.5   & 5.0E-11 & 7 \\ \hline
		$K_1$ & Dirichlet & 0.5   & 5.4E-14 & 8  \\ \hline
	\end{tabular}%
	\label{tab:test_on_convergence}%
\end{table}%
\subsection{Numerical results on performance of kernel with singularities}\label{subsec:numerical_results_sigularity}
This subsection includes the similar testing results as in Subsection \ref{subsec:numerical_results_continous} for FIDEs with kernel $K_1$ that bears singularities when $\gamma<0$. 

Table \ref{tab:Fredholm_en_gamma} shows the impact of degree of singularity of $k(x,y)$. Similar performance is observed for all covered $\gamma$, consistent to what has been observed in Subsection \ref{subsec:numerical_results_continous}.
\begin{table}[htbp]
	\caption{Impact of kernel's smoothness with $\mu=1,\theta=3\pi/2, M=2^q$}
	\begin{tabular}{lrrrr}
		kernel & type & \multicolumn{1}{l}{$\gamma$} & \multicolumn{1}{l}{$max_{err}$} & \multicolumn{1}{l}{$q$} \\ \hline\hline
		$K_1$ & Dirichlet & -0.9   & 3.0E-08 & 7 \\\hline
		$K_1$ & Dirichlet & -0.5   & 2.7E-08 & 7 \\\hline
		$K_1$ & Dirichlet & 0     & 3.0E-08 & 7 \\\hline
		$K_1$ & Dirichlet & 0.5  & 3.5E-08 & 7 \\\hline
		$K_1$ & Dirichlet & 1.5  & 4.8E-08 & 7 \\\hline
		$K_1$ & Dirichlet & 2    & 5.8E-08 & 7 \\ \hline
	\end{tabular}%
	\label{tab:Fredholm_en_gamma}%
\end{table}%
Table \ref{tab:Fredholm_en_convergence} shows convergence property and high accuracy is achieved for $q\ge 7$ with Dirichlet condition. As expected, the performance is improved as $q$ increases. 
\begin{table}[htbp]
	\caption{Convergence Test with $f(x)=\cos (\frac{3\pi x}2)$}
	\begin{tabular}{lrrrr}
		kernel & type & \multicolumn{1}{l}{$\gamma$} & \multicolumn{1}{l}{$max_{err}$} & \multicolumn{1}{l}{$q$} \\ \hline\hline
		$K_1$ & Dirichlet & -0.5   & 4.9E-03 & 4 \\\hline
		$K_1$ & Dirichlet & -0.5   & 9.9E-05 & 5 \\\hline
		$K_1$ & Dirichlet & -0.5   & 4.6E-07 & 6 \\\hline
		$K_1$ & Dirichlet & -0.5   & 2.7E-08 & 7 \\\hline
		$K_1$ & Dirichlet & -0.5   & 1.7E-09 & 8 \\\hline
		$K_1$ & Dirichlet & -0.5   & 1.1E-10 & 9 \\ \hline
	\end{tabular}%
	\label{tab:Fredholm_en_convergence}%
\end{table}%
Table \ref{tab:Fredholm_en_boundary} shows the performance of different boundary conditions with target function $f(x)=\cos(3\pi x/2)$ for FIDE with kernel $K_1, \gamma=-0.5$. The performance is decent for all cases and reaches the best accuracy with Dirichlet condition, consistent to what is observed in Subsection \ref{subsec:numerical_results_continous}.
\begin{table}[htbp]
	\caption{The performance of covered boundary conditions with   $f(x)=\cos (\frac{3\pi x}2)$}
	\begin{tabular}{llrrr}
		kernel & type & \multicolumn{1}{l}{$\gamma$} & \multicolumn{1}{l}{$max_{err}$} & \multicolumn{1}{l}{$q$} \\ \hline\hline
		$K_1$ & Neumann & -0.5   & 1.6E-07 & 7 \\\hline
		$K_1$ & Dirichlet & -0.5   & 2.7E-08 & 7 \\\hline
		$K_1$ & $Mix_1$ & -0.5   & 9.7E-08 & 7 \\\hline
		$K_1$ & $Mix_2$ & -0.5   & 6.4E-08 & 7 \\ \hline
	\end{tabular}%
	\label{tab:Fredholm_en_boundary}%
\end{table}%
Table \ref{tab:Fredholm_en_boundary} provides the performance of FIDE with kernel $K_1, \gamma=-0.5$ with covered target functions.
%and Figure \ref{fig:multiple_idef_en} compares numerical solutions with target functions over the range $[0,b]$. 
High accuracy is reached for all cases. 
\begin{table}[htbp]g
	\caption{Performance with covered target functions}
	\begin{tabular}{llrrrr}
			kernel & type & \multicolumn{1}{l}{$\gamma$} &  target function & \multicolumn{1}{l}{$max_{err}$}  & \multicolumn{1}{l}{$q$}  \\ \hline\hline
		$K_1$ & Dirichlet & -0.5  & $cos(3\pi x/2)$    &  2.7e-08  & 7 \\\hline
		$K_1$ & Dirichlet & -0.5  & $cos(\pi x/2)$    & 9.7e-09  & 7 \\\hline
		$K_1$ & Dirichlet & -0.5  & $x^2$    & 6.9e-09  & 7 \\\hline
		$K_1$ & Dirichlet & -0.5  & $exp(x)$    & 9.3e-10 & 7 \\  \hline
	\end{tabular}%
	\label{tab:Fredholm_en_multiple}%
\end{table}%
\section{Summary}\label{sec:summary}
In this paper, we enhance the trigonometric interpolation algorithm in \cite{zou_tri} to approximate non-periodic functions in $2$ dimensional space.   Performance of the enhanced Algorithm \ref{alg:dim2} is expected to be similar as in $1$-dim case and achieves accuracy aligned with smoothness of underlying function, which has been confirmed by numerical tests conducted in Section \ref{sec:trig_dim2}. 

Algorithm \ref{alg:dim2} has been applied to develop an approximation method, Algorithm \ref{Alg:FIDE}, for the solution of a second order linear Fredholm integro-differential equation (FIDE).  The new method converts FIDE to a linear algebraic system as for the treatment on second order linear ODE in \cite{zou_tri_IV} and it bears several advantages as mentioned in Section \ref{sec:intro}. The tests on Algorithm \ref{alg:dim2} in Section \ref{sec:numerical_examples} include various boundary conditions with both continuous and integrable kernels.  Decent performance is observed across all covered scenarios with moderate size of grid points.


\begin{thebibliography}{2}

\bibitem{ide_23}
K. Maleknejad, and Y. Mahmoidi,
\textit{Numerical solution of linear Fredholm integral equation by using hybrid Taylor and block-pulse functions}, Appl. Math. Comput., 149 (2004), pp. 799-806.


\bibitem{ide_27}
M.T. Rashed,
\textit{Numerical solution of functional differential, integral and integro-differential equations}, Appl. Numer. Math., 156 (2004), pp 485-492. 

\bibitem{ide_30}
Y. Ren, B. Zhang, H. Qiao,
\textit{A simple Taylor-series expansion method for a class of second kind integral equations.}, J. Comp. Appl. Math. 110 (1999), pp 15--24. 

\bibitem{ide_36}
W. Wang. C. Lin.
\textit{A new algorithm for integral of trigonometric functions with mechanization}, Appl. Math. Comput., 164 (2005), pp 71-82.

\bibitem{ide_38}
S. Yal¸cınba¸s, M. Sezer,
\textit{The approximate solution of high-order linear Volterra-Fredholm integro-differential equations in terms of Taylor polynomials}, Appl. Math. Comput., 164, 2000, pp 291-308

\bibitem{ide_10}
G. Ebadi, M.Y. Rahimi, S. Shahmorad,
\textit{Numerical solution of the nonlinear Volterra integro-differential equations by the Tau method.}, {Appl. Math. Comput.}, 188 (2007), pp {1580-1586},

\bibitem{ide_28}
M. Razzaghi, S. Yousefi,
\textit{Legendre wavelets method for the nonlinear Volterra Fredholm integral equations.}, {Math. Comput. Simul.}, 70 (2005), pp {1-8}

\bibitem{ide_1}
G. Akram, S. Siddiqi,
\textit{Solution of sixth order boundary value problems using non-polynomial spline technique.}, {Appl. Math. Comput.}, 181 (2006), pp 708-720,


\bibitem{ide_4}
A. Akyuz-Dascioglu,
\textit{A Chebyshev polynomial approach for linear Fredholm-Volterra integro-differential equations in the most general form.}, {Appl. Math. Comput.}, 181 (2006), pp {103-112},


\bibitem{ide_6}
N.K. Basu,
\textit{A Chebyshev series method for the numerical solution of Fredholm integral equations with associated eigen-value problem.}, {SIAM J. Numer. Anal.}, 10 (1971), pp 57-68


% 85 (1),
\bibitem{ode_desa}
M. Dehghan and A. Saadatmandi,
\textit{Chebyshev finite difference method for Fredholm integro-differential equation}, {Int. J. Comput. Math.}, 85 (2008), pp 123-130


\bibitem{zou_tri}
X. Zou
\textit{On Trigonometric Approximation and Its Applications}, \qquad\qquad\qquad\qquad\qquad\qquad\qquad\qquad https://arxiv.org/pdf/2505.02330. May, 2025

\bibitem{zou_tri_III}
X. Zou
\textit{Trigonometric Interpolation Based Optimization for Solving Second Order Non-Linear ODE  with Linear Boundary Conditions}, \qquad\qquad\qquad\qquad\qquad\qquad\qquad\qquad
https://arxiv.org/pdf/2504.19280, Apr. 2025


\bibitem{zou_tri_IV}
X. Zou
\textit{Trigonometric Interpolation Based Approach for Second Order ODE with Mixed Boundary Conditions}, 
https://arxiv.org/pdf/2505.07183,  May, 2025


\bibitem{zou_tri_VIDE}
X. Zou,
\textit{Trigonometric Interpolation Based Approach for Second Order Volterra integro-differential equations}, to be appear


\end{thebibliography}
\end{document}